\documentclass{amsart}
\usepackage{amsmath, amssymb, graphicx, amsthm, relsize}
\usepackage{hyperref}

\usepackage[sorted, msc-links, backrefs]{amsrefs}

\usepackage[margin=1.3in]{geometry}

\newtheorem {Theorem}{Theorem}[section]

\newtheorem {Remark}[Theorem]{Remark}

\newtheorem {Lemma}[Theorem]{Lemma}

\newcommand{\beq}{\begin{equation}}
\newcommand{\eeq}{\end{equation}}

\def\bC{{\mathbb C}}
\def\bF{{\mathbb F}}
\def\c{{\bf c}}
\def\e{{\bf e}}
\def\bQ{{\mathbb Q}}
\def\Q{{\mathbb Q}} \def\dd{{\bf d}}
\def\bZ{{\mathbb Z}} 
\def\H{{\mathcal H}} 
\def\Hred{\H_{red}}
\def\bfC{{\mathcal C}} 
\def\A{{\bf A}}
\def\la{{\nu}} 
\def\bp{{\bf p}}
\def\OR{{\mathcal U}^{(r)}} 
\def\G{{\Sigma}}
\def\bN{{\mathbb N}}
\def\Qab{{\bQ_{\rm ab}}}
\def\bR{{\mathbb R}}
\def\aa{{\alpha}}
\def\t{\tilde}
\def\tg{\t\g}
\def\p{\pi_1 (X,a)}
\def\f{f:R \to S}
\def\E{\bf E}
\def\Est{\mathbb D ^*}
\def\F{{\mathcal F}}
\def\P{{\mathbb P}^1} \def\xs{{\mathbb o}}
\def\M{{\mathcal M}(Y)}
\def\Ma{{\mathcal M}(Y^{(1)})}
\def\Mb{{\mathcal M}(Y^{(2)})}
\def\eps{\varepsilon}
\def\D{\Delta}
\def\ur{{\mathcal U}_r}
\def\Dh{\hat{\D}}

\def\fge{{$L/k(x)$\ }} 
\def\bC{{\mathbb C}} 
\def\t{{\tau}}
\def\k{\kappa} 
\def\L{\Lambda} 
\def\C{{\bf C}} 
\def\D{{\bf D}}
\def\bF{{\mathbb F}}
\def\bt{{\bf t}} 
\def\bfC{{\bf C}}
\def\HC{{\mathcal H}_r(\C)} 
\def\S{{\cal S}}
\def\SA{{\mathcal S}_r^{(\A)}(G)} 
\def\orr{{\cal O}(r+1)}
\def\TG{{\mathcal T}_r(G)} 
\def\oor{{\cal O}_r} 
\def\T{{\cal T}}
\def\K{{\mathcal K}} 
\def\k{\kappa}
\def\A{{\bf A}} 
\def\UR{{\cal U}^{(r)}}
\def\Ni{{\mbox{Ni}}} 
\def\Vtt{(\tilde V)\tilde{}}
\def\Vht{(\hat V)\tilde{}}
\def\d{{\delta }} 
\def\a{{\alpha }} 
\def\b{{\beta }} 
\def\g{{\gamma }}
\def\M{{\mathcal M}} 
\def\FF{{\cal F}} 
\def\s{{\sigma }}
\def\H{{\mathcal H}} 
\def\L{{\mathcal L}} 
\def\NN{{\cal N}}
\def\HHR{{\mathcal H}^{(r)}} 
\def\HHr{{\cal H}_{r}}
\def\U{{\mathcal U}} 
\def\UR{{\cal U}^{(r)}} 
\def\Ur{{\cal U}_{r}}
\def\Aut{\mbox {Aut}}

\def\endofproof{\hfill\vrule height6pt width5pt depth0pt}
\def\Proof{{\bf Proof :}~}

\hypersetup{
  colorlinks   = true,	
  urlcolor     = blue,
  linkcolor    = blue,
  citecolor   = red
}

\title{The locus of curves with prescribed automorphism group}

\begin{document}

\maketitle

\centerline{April 24, 2002\footnote{This papers is exactly the same version of  its original version posted on 2002.  The only difference is that this version has links to all the XIX-century references, which were not available at the time of publication.}}
\bigskip

\centerline{K. Magaard\footnote{Partially supported by NSA grant MDA904-01-1-0037}, T.Shaska, S.
Shpectorov\footnote{Partially supported by NSA grant   MDA904-01-1-0023}, and H.
V\"olklein\footnote{Partially supported by NSF grant DMS-9970357 }}

\bigskip

\begin{center}Wayne State University, University of California at Irvine,\\
Bowling Green State University, and University of Florida
\end{center}

\begin{abstract} Let $G$ be a finite group, and $g\ge2$. We study the locus of genus
$g$ curves that admit a $G$-action of given type, and inclusions between such loci. We use this to study the
locus of genus $g$ curves with prescribed automorphism group $G$. We completely classify these loci for $g=3$
(including equations for the corresponding curves), and for $g\le10$ we classify those loci corresponding to
``large'' $G$.
\end{abstract}

\section{Introduction}
There is a vast literature on automorphism groups of compact Riemann surfaces, beginning in the 19th century
with Schwartz, Klein, Hurwitz, Wiman and others. However, most of the literature is quite recent. In the
first part of the paper, we survey the main results.

By covering space theory, a finite group $G$ acts (faithfully) on a genus $g$ curve if and only if it has a
genus $g$ generating system (see section \ref{fullsection} below). Using this purely group-theoretic
condition, Breuer \cite{Br} classified all groups that act on a curve of genus $\le 48$. This was a major
computational effort using the computer algebra system GAP \cite{GAP}. It greatly improved on several papers
dealing with small genus, by various authors (see the references in part I).

Of course, for each group in Breuer's list, all subgroups are also in the list. This raises the question how
to pick out those groups that occur as the {\bf full automorphism group} of a genus $g$ curve. This question
is answered in Part II of the paper. Let $\M_g$ be the moduli space of genus $g$ curves. We study the locus
$\L$ in $\M_g$ of curves admitting a $G$-action of given ramification type (resp., signature). All components
of $\L$ have the same dimension which depends only on the signature of the $G$-action. Restricting the action
to a subgroup $H$ of $G$ yields an inclusion of $\L$ into the corresponding locus $\L'$ for the action of
$H$. If $\dim \L<\dim L'$ then for "most" curves in $\L'$, the $H$-action does not extend to a $G$-action of
the type defining $\L$. Thus one is led to classify the pairs $(\L,\L')$ with $\dim \L=\dim\L'$. This is done
in Lemma \ref{Lemma2} below: It turns out that such pairs exist only in very restricted cases, in particular
only if $\dim\L\le3$. From that we derive a necessary and sufficient condition for a group to occur as the
full automorphism group of a genus $g$-curve (Theorem \ref{thm2}).

After Part II was written, we found references to papers of Singerman \cite{Si} and Ries \cite{Ri} which
contain similar results. Their method is different, of analytic nature, using Teichm\"uller theory and
Fuchsian groups. Our approach is based on algebraic geometry, using the algebraic structure of Hurwitz spaces
and the moduli spaces $\M_g$; therefore, it can be used to obtain information on fields of definition. This
aspect may be studied in later work.

In Part III we apply the above criterion to the data compiled by Breuer (available from his website, see 1.5
below). Our first application is in the case of genus 3, which is already quite rich and shows the power of
our group-theoretic method. We obtain the full list of automorphism groups in genus 3 plus equations for the
corresponding curves. This result is scattered over several long papers by Kuribayashi and various co-authors
\cite{KuKo1}, \cite{KuKo2}, \cite{KK1}, \cite{KK2}, some of which contain errors (see section 6 below) and
some are not available in most libraries. None of them seems to give a complete account of this basic result.
Here we show how to derive it quickly from the group-theoretic data.

Further, we obtain the list of ``large'' groups $\Aut(X_g)$ (see 1.3 below) up to $g=10$ and the dimension and number of components of the corresponding loci in $\M_g$, see Table 4.

The locus of curves in $\M_g$ with fixed automorphism group consists of finitely many components; to
determine their number requires tools that will be developed in later work (mapping class group action on
generating systems). In the special case of components corresponding to curves of orbit genus $0$ (i.e., the
quotient by the automorphism group of the curve has genus $0$), the mapping class group action is just the
braid group action studied in \cite{FV}, \cite{V}, \cite{Buch}; and those components correspond to braid
orbits on generating systems, which can be computed with our GAP package BRAID \cite{MSV}. A curve with large
automorphism group always has orbit genus $0$, so we were able to use the BRAID package to compute the
information in Table 4. The details of this computation will be given in later work (in a more general
situation).

\medskip

\noindent  \textbf{Notation:} We will use the term ``curve'' to mean ``compact Riemann surface''; and $X=X_g$ denotes a
curve of genus $g\ge2$. Further, $\Aut(X)$ denotes the group of analytic (equivalently, algebraic)
automorphisms of $X$.

\part{Survey of known results}
\subsection{The beginnings in the 19th century}
The group $\Aut(X)$ acts on the finite set of Weierstra\ss\ points of $X$. This action is faithful unless $X$
is hyperelliptic, in which case its kernel is the group of order 2 containing the hyperelliptic involution of
$X$. Thus in any case, $\Aut (X)$ is a finite group. This was first proved by Schwartz.

F. Klein \cite{Kl1} (1879) studied the genus 3 projective curve
\begin{equation}
x^3y\ +\ y^3z\ +\ z^3x\ \ =\ \ 0 \label{Klein}
\end{equation}
which is now known as Klein's quartic. He found all its automorphisms: They form the simple group of order
168. Further, he showed that there is a degree 168 covering from this curve to $\P$ branched over three
points with ramification indices 2,3,7. See \cite{Le} for a detailed history on Klein's quartic (including an
English translation of Klein's original article) as well as an account of further developments influenced by
it.

Poincar\'e introduced Fuchsian groups in 1882. They are basic for the analytic theory of Riemann surfaces and
have been used heavily in the study of $\Aut(X)$ (see \cite{Br}, \cite{Ri} for example). In this paper we
take a different (algebraic) point of view.

The next milestone was Hurwitz's seminal paper \cite{Hu} in 1893, where he discovered what is now called the
Riemann-Hurwitz formula. From this he derived that
$| \Aut (X_g)| \ \ \leq \ \ 84 \, \, (g-1)$, now known as the Hurwitz
bound. This bound is attained by Klein's quartic, for example. Curves which attain this bound are called {\bf
Hurwitz curves} and their automorphism groups {\bf Hurwitz groups}. Hurwitz proved that a finite group is a
Hurwitz group if and only if it has generators $a,b,c$ of orders $2,3,7$, respectively, with $abc=1$.

\subsection{Hurwitz groups}
Klein's quartic is the only Hurwitz curve of genus $g\le3$. Fricke showed that the next Hurwitz group occurs
for $g=7$ and has order 504. Its group is $SL(2,8)$, and an equation for it was computed by Macbeath
\cite{Mb2} in 1965. Klein's quartic and Macbeath's curve are the only Hurwitz curves whose equations are
known. Further Hurwitz curves occur for $g=14$ and $g=17$ (and for no other values of $g\le19$).

There are a lot of papers by group-theoretists on Hurwitz groups, surveyed by Conder \cite{Co}. It follows
from Hurwitz's presentation that a Hurwitz group is {\bf perfect}. Thus every quotient is again a Hurwitz
group, and if such a quotient is minimal then it is a non-abelian simple group. Several infinite series of
simple Hurwitz groups have been found by Conder, Malle, Kuribayashi, Zalessky, Zimmermann and others. In
2001, Wilson \cite{Wl} showed the monster is a Hurwitz group.

\subsection{Large automorphism groups.}
For a fixed $g\geq 2$ denote by $N(g)$ the maximum of the $|\Aut(X_g)|$. Accola \cite{Ac1} and Maclachlan
\cite{Mc1} independently show that $N(g) \geq 8 (g+1)$ and this bound is sharp for infinitely many $g$'s. If
$g$ is divisible by 3 then $N(g) \geq 8 (g+3)$.

The following terminology is rather standard. We say $G\le\Aut(X_g)$ is a {\bf large automorphism group} in
genus $g$ if  $|G|\ \ >\ \ 4 (g-1)$.
Then the quotient of $X_g$ by $G$ is a curve of genus $0$, and the number of points of this quotient ramified
in $X_g$ is 3 or 4 (see \cite{Br}, Lemma 3.18, or \cite{FK}, pages 258-260). Singerman \cite{S1} (1974) shows
that Riemann surfaces with large cyclic, Abelian, or Hurwitz groups are symmetric (admit an involution).
Kulkarni \cite{Ku}(1997) classifies Riemann surfaces admitting large cyclic automorphism groups and works out
full automorphism groups of these surfaces. Matsuno \cite{Mt}(1999) investigates the Galois covering of the
projective line from compact Riemann surfaces with large automorphism groups.

\subsection{Cyclic groups as automorphism groups.}
Let $t $ be the order of an automorphism of $X_g$. Hurwitz \cite{Hu} showed $ t\leq 10 (g-1)$. In 1895, Wiman
improved this bound to be $t\le 2(2 g+1)$ and showed this is best possible. If $t$ is a prime then $t \leq 2
g +1$. Homma \cite{Ho} (1980) shows that this bound is achieved if and only if the curve is birationally
equivalent to
\[
y^{m-s} (y-1)^s = x^q, \quad for \quad 1 \leq s < m \leq g_x+1
\]

\subsection{ The canonical representation of $\Aut(X)$}
Each group $G\le \Aut(X_g)$ acts faithfully on the $g$-dimensional vector space of holomorphic differential
forms on $X$. The resulting finite subgroups of $GL(g,\bC)$ satisfy a number of special conditions, studied
by I. Kuribayashi, A. Kuribayashi, Kimura, Ohmori, Kobayashi, Hayakawa and others. A. Kuribayashi and others
(see \cite{KuKi1} and the references in that paper) have completely listed these finite linear groups up to genus
5.

Using computational group theory (more precisely, the computer algebra system GAP \cite{GAP}), Breuer
\cite{Br} extended this up to genus 48: Computing all groups $G\le\Aut(X_g)$ together with their character on
the space of holomorphic differential forms. He avoids ad-hoc methods by using the GAP-library of ``small
groups'' by Besche and Eick \cite{BE}. His results are collected in a database which is available from his
website at 
\href{www.math.rwth-aachen.de/~Thomas.Breuer/}{www.math.rwth-aachen.de/~Thomas.Breuer/}. 
The GAP library
of small groups is available at 
\href{http://www-public.tu-bs.de:8080/\~{}beick/so.html}{http://www-public.tu-bs.de:8080/\~{}beick/so.html}. 

In our tables a group $G$ is identified by its Group
ID from the Small Groups Library. A Group ID is a pair $(n,m)$ where $n$ is the order of $G$ and $m$ the
number of $G$ among the groups of order $n$.

\subsection{Full automorphism groups.}
For each group in Breuer's list, all subgroups are also in the list. This raises the question how to pick out
those groups that occur as the {\bf full automorphism group} of a genus $g$ curve. Answering this question
requires a moduli space argument: Either analytically, using Fuchsian groups and Teichm\"uller space, as done
by Singerman \cite{Si} and Ries \cite{Ri}, or algebraically, as done in Part II of this paper. Using these
general results, one can compute explicitly all groups $\Aut(X_g)$ for small $g$. We give the complete list
for $g=3$ in Part III, and the list of large groups $\Aut(X_g)$ up to $g=10$. Further results will be given
in later work, together with additional information.

For $g>3$, there are no lists of full automorphism groups in the literature. The case $g=2$ is easy, but
already for $g=3$ there is no complete account (see the remarks in the introduction).

\part{Subgroups of $\Aut(X_g)$ not occurring as full automorphism groups in genus $g$}

\section{Ramification type and signature of a $G$-curve}
Fix an integer $g\ge2$ and a finite group $G$. Let $C_1,...,C_r$ be conjugacy classes $\ne\{1\}$ of $G$. Let
$\C=(C_1,...,C_r)$, viewed as unordered tuple, repetitions are allowed. We allow $r$ to be zero, in which
case $\C$ is empty.

We will use the term ``curve" to mean ``compact Riemann surface". Consider pairs $(X,\mu)$, where $X$ is a
curve and $\mu: G\to\Aut(X)$ is an injective homomorphism. Mostly we will suppress $\mu$ and just say $X$ is
a curve with $G$-action, or a $G$-curve, for short. Two $G$-curves $X$ and $X'$ are called equivalent if
there is a $G$-equivariant isomorphism $X\to X'$.

We say a $G$-curve $X$ is {\bf of ramification type} $(g,G,\C)$ if the following holds: Firstly, $g$ is the
genus of $X$. Secondly, the points of the quotient $X/G$ that are ramified in the cover $X\to X/G$ can be
labelled as $p_1,...,p_r$ such that $C_i$ is the conjugacy class in $G$ of distinguished inertia group
generators over $p_i$ (for $i=1,...,r$). (Distinguished inertia group generator means the generator that acts
in the tangent space as multiplication by $\exp(2\pi\sqrt{-1}/e)$, where $e$ is the ramification index). For
short, we will just say $X$ is of type $(g,G,\C)$.

If $X$ is a $G$-curve of type $(g,G,\C)$ then the genus $g_0$ of $X/G$ is given by the Riemann-Hurwitz
formula
\[
\frac{2\ (g-1)}{|G|}\ \ \ = \ \ \ 2\ (g_0-1)\ +\ \sum_{i=1}^r\ \left(1-\frac{1}{c_i} \right) \leqno{(1)},
\]
 where $c_i$ is
the order of the elements in $C_i$.

Note that $g_0$ (the {\bf orbit genus}) depends only on $g$, $|G|$ and the {\bf signature}
$\c=(c_1,\ldots,c_r)$ of the $G$-curve $X$.

\section{Hurwitz spaces and moduli of curves}\label{section2}

Define $\H=\H(g,G,\C)$ to be the set of equivalence classes of $G$-curves of type $(g,G,\C)$. By covering
space theory (or the theory of Fuchsian groups), $\H$ is non-empty if and only if $G$ can be generated by
elements $\a_1,\b_1,...,\a_{g_0},\b_{g_0},\g_1,...,\g_r$ with $\g_i\in C_i$ and
\[
\prod_{j}\ [\a_j,\b_j]\ \ \prod_{i}\ \g_i\ \ \ \ = \ \ \ 1 \leqno{(2)}
\]
 Here $[\a,\b]=\ \a^{-1}\b^{-1}\a\b$.

Let $\M_{g}$ be the moduli space of genus $g$ curves, and $\M_{g_0,r}$ the moduli space of genus $g_0$ curves
with $r$ distinct marked points, where we view the marked points as unordered (contrary to usual procedure).
Consider the map
\[\Phi:\ \H\ \to \ \M_{g}\]
forgetting the $G$-action, and the map
\[
\Psi:\ \H\ \to \ \M_{g_0,r} 
\]
mapping (the class of) a $G$-curve $X$ to the class of the quotient curve $X/G$ together with the (unordered)
set of branch points $p_1,...,p_r$. If $\H\ne\emptyset$ then $\Psi$ is surjective and has finite fibers, by
covering space theory. Also $\Phi$ has finite fibers, since the automorphism group of a curve of genus $\ge2$
is finite.

By \cite{B}, the set $\H$ carries a structure of quasi-projective variety (over $\bC$) such that the maps
$\Phi$ and $\Psi$ are finite morphisms. If $\H\ne\emptyset$ then all components of $\H$ map surjectively to
$\M_{g_0,r}$ (through a finite map), hence they all have the same dimension
\[
\d(g,G,\C):= \ \ \dim\ \M_{g_0,r} \ \ = \ \ 3g_0-3+r
\]
(This is correct because $r\ge1$ if $g_0=1$ and $r\ge3$ if $g_0=0$; the latter holds because $g\ge2$). Since
also $\Phi$ is a finite map, we get

\begin{Lemma} \label{Lemma1} Let $\M(g,G,\C)$ denote the image of $\Phi$,
i.e., the locus of genus $g$ curves admitting a $G$-action of type $(g,G,\C)$. If this locus is non-empty
then each of its components has dimension $\d(g,G,\C)$.
\end{Lemma}

Note that $\d(g,G,\C)$ depends only on $g$, $|G|$ and the signature, so we write $\d(g,G,\c):=\d(g,G,\C)$.
The Lemma continues to hold if we replace $\C$ by $\c$.

\section{Restriction to a subgroup} \label{restr}

Let $H$ be a subgroup of $G$. Then each $G$-curve can be viewed as an $H$-curve by restriction of action. Let
$X$ be a $G$-curve of type $(g,G,\C)$. Then the resulting $H$-curve is of type $(g,H,\D)$, where $\D$ is
obtained as follows: Choose $\g_i\in C_i$ and let $\s_{i,1}, \s_{i,2},...$ be a set of representatives for
the double cosets $<\g_i>\s H$ in $G$. Let $m_{ij}$ be the smallest integer $\ge1$ such that the element
$\s_{ij}^{-1}\g_{i}^{m_{ij}}\s_{ij}$ lies in $H$, and let $D_{ij}$ be the conjugacy class of this element in
$H$. Then $\D$ is the tuple consisting of all $D_{ij}$ with $D_{ij}\ne\{1\}$. (More precisely, the tuple $\D$
is indexed by the set of possible pairs $(i,j)$, and its $(i,j)$-entry is $D_{ij}$ ). The definition of $\D$
does not depend on the choice of the $\g_i$ and $\s_{ij}$. Note that the signature of the $H$-curve depends
on the type of the $G$-curve, not only on its signature.

We have $$ \M(g,G,\C) \ \ \subset \ \ \M(g,H,\D)$$ hence their dimensions satisfy
\[
\d(g,G,\C)\ \ \ \le\ \ \ \d(g,H,\D)
\]
If this is a strict inequality then the complement of the closure of $ \M(g,G,\C)$ in $\M(g,H,\D)$ is open
and dense; then in particular, it is not true that every $H$-curve of type $(g,H,\D)$ is the restriction of a
$G$-curve of type $(g,G,\C)$.

\begin{Lemma} \label{Lemma2} Suppose $H$ is a subgroup of $G$ of index
$n>1$. We let $\bar G$ denote the permutation group induced by $G$ on the coset space $G/H$. Let $X$ be a
$G$-curve of type $(g,G,\C)$, where we label the $C_i$ such that $c_1\le \ldots \le c_r$. Let
$\c=(c_1,\ldots,c_r)$. The $H$-curve obtained by restriction of action is of type $(g,H,\D)$, where $\D$ is
defined above; let $(g,H,\dd)$ be its signature. Let $s$ (resp., $r$) be the length of $\D$ (resp., $\C$).
Let $h_0$ (resp., $g_0$) be the genus of $X/H$ (resp., $X/G$). If $\d:=\d(g,G,\C)=\d(g,H,\D)$ then $g_0=0$,
and one of the following holds:
\begin{enumerate}
\item[(I)]\ \ $\d=3$: \ Then $n=2$, $h_0=2$, $s=0$, $r=6$, all $c_i=2$ and $C_i\not \subset H$.

\item[(II)]\ \ $\d=2$:\ \ Then $n=2$, $h_0=1$, $s=2$, $r=5$, $C_i\not\subset H$ and $c_i=2$ for
$i=1,...,4$ and\ $\dd=(c_5,c_5)$.
\item[(III)]\ \ $\d=1$:

\subitem (a) \ $n=2$, $h_0=1$, $s=1$, $r=4$, \ $C_i\not\subset H$ for all $i$, \ $c_1=c_2=c_3=2$ and \
$\dd=(c_4/2)$.

\subitem (b) \ $n=2$, $h_0=0$, $s=4=r$, $c_1=c_2=2$, $C_3,C_4\subset H$, $\dd=(c_3,c_3,c_4,c_4)$.

\subitem (c) \ $h_0=0$, $s=4=r$, $c_1=c_2=c_3=2$, $C_4\subset H$, and $H$ is normal in $G$ with $G/H$ a Klein
4-group; 
and $\dd=(c_4,c_4,c_4,c_4)$.

\item[(IV)]\ \ \ $\d=0$: \ Then $h_0=0$, $s=3=r$ and

\subitem (a) \ $n=2$, $h_0=0$, $r=3=s$, $c_1=2$, $c_i>2$ is even and $C_j\subset H$ where $\{i,j\}=\{2,3\}$;
\subitem \hskip .3in and $\dd=(c_i/2,c_j,c_j)$.

\subitem (b)\ $n=3$, $\bar G\cong \bZ/3$, one $C_i\subset H$ and the other two have $c_j=3$; here
$\dd=(c_i,c_i,c_i)$. \subitem (c) \ $n=3$, $\bar G\cong S_3$, $c_1=2$, $c_2=3$, $c_3>2$ is even,
$\dd=(2,c_3/2,c_3)$. \subitem (d) \ $n=4$, $\bar G\cong A_4$, $c_1=2$, $c_2=3$, $c_3>3$ is divisible by 3,
$\dd=(3,c_3/3,c_3)$. \subitem (e) \ $n=6$, $\bar G\cong \mbox{PGL}_2(5)$, $\c=(2,4,5)$, $\dd=(4,4,5)$,
$g\equiv 1$ mod 3. \subitem (f) \ $n=8$, $\bar G\cong \mbox{PSL}_2(7)$, $\c=(2,3,7)$, $\dd=(3,3,7)$, $g\equiv
1$ mod 2. \subitem (g) \ $n=9$, $\bar G\cong \mbox{PSL}_2(8)$, $\c=(2,3,7)$, $\dd=(2,7,7)$, $g\equiv 1$ mod
6. \subitem (h) \ $n=10$, $\bar G\cong \mbox{PGL}_2(9)$, $\c=(2,3,8)$, $\dd=(3,8,8)$, $g\equiv 1$ mod 15.
$g\equiv 1$
\subitem (i) \ There are 6 further cases when $H$ is not maximal in $G$, and $n=\ 4,6,6,12,12,24$,
respectively.
\end{enumerate}
\end{Lemma}

\noindent The cases in (IV) are taken from \cite{M}, Table 2 (c.f. 5b below).

\smallskip\noindent
\Proof We denote the hypothesis $\d(g,G,\C)=\d(g,H,\D)$ by (Hyp). Let $r_1$ be the number of points of $X/G$
that are ramified in $X/H$. Applying Riemann-Hurwitz to the cover $X/H\to X/G$ of degree $n$ yields
$$ 2\ (h_0-1) \ \ \ge \ \ 2n\ (g_0-1) \ +\ r_1 \leqno{(3)}$$
Let $s$ be the number of points of $X/H$ that are ramified in $X$, and $r_2:=r-r_1$. The $r_2$ points of
$X/G$ that are ramified in $X$, but not in $X/H$, have a total of $nr_2$ distinct pre-images in $X/H$. The
latter points are ramified in $X$, hence $s\ge nr_2$.

\smallskip\noindent
{\bf Case 1:} \ $g_0>1:$  \;  Then (Hyp) gives
$$ 3g_0\ +\ r_1 \ +\ r_2 \ \ = \ \ 3h_0\ +\ s \ \ \ge \ \ 3h_0\ +\ r_2$$
hence $ 3g_0\ +\ r_1 \ \ \ge \ \ 3h_0$.  With (3) this gives
$$3g_0\ +\ 2\ (h_0-1) \ -\ 2n\ (g_0-1)\ \ \ge \ \ 3h_0$$
$$ g_0\ (3-2n)\ +\ 2n \ \ \ge \ \ h_0 \ +\ 2 \ \ \ge \ \ g_0 \ +\ 2 $$
$$g_0\ \ <\ \ \frac{n}{n-1}\ \ \le \ \ 2$$
a contradiction.

\smallskip\noindent
{\bf Case 2:} \ $g_0=1:$  \;  If $r=0$ then by (1) we have $g=1$, contrary to assumption. Thus $r>0$. Now
(3) gives $ 2\ (h_0-1) \ge r_1$. If $h_0>1$ then (Hyp) gives $r_1\ge\ 3(h_0-1)$ (as in Case 1), a
contradiction. If $h_0=1$ then $r_1=0$, hence (Hyp) gives $r_2=s\ge nr_2$. Contradiction, since $r_2=r>0$.

\medskip
We have now proved $g_0=0$. Thus if $h_0>0$ then $r_1\ge3$. Further, we can choose $\g_i\in C_i$ with
$\g_1\cdots\g_r=1$ and $G=<\g_1,...,\g_r>$ (by (2)). Then $c_i$ is the order of $\g_i$. Let $p_1,...,p_r$
(resp., $q_1,...,q_{s}$) be the points of $X/G$ (resp., $X/H$) that are ramified in $X$. We choose the
labelling such that $p_i$ corresponds to $C_i$ and
$$c_1\ \ \le \ \ldots \ \le \ \ c_r\leqno{(4)}$$
The orbits of $\g_i$ on $G/H$ correspond to the points of $X/H$ over $p_i$ (where the length of the orbit
equals the ramification index of the corresponding point over $p_i$). Let $\mbox{Ind}(C_i)$
(=$\mbox{Ind}(\g_i)$ ) be $n$ minus the number of those orbits. Those orbits that correspond to points
$\not\in\{q_1,...,q_{s}\}$ have length $c_i\ge2$. Thus the total number of orbits of all $\g_i$ is $\le
s+(nr-s)/2$. Hence Riemann-Hurwitz for $X/H\to X/G$ yields
$$2\ (h_0-1+n) \ \ \ =\ \ \ \sum_{i=1}^{r}\ \mbox{Ind}(C_i) \ \ \ \ge\ \ \   nr\ - \ \frac{nr+s}{2} \ \ \ =\ \ \ \frac{nr-s}{2}\leqno{(5)}$$

\medskip\noindent
{\bf Case 3:} \ $h_0\ge2:$ \;  Here (Hyp) gives $r=3h_0+s$. Plugging this into (5) gives
$$2\ (h_0-1+n) \ \ \ \ge\ \ \ \frac{n\ (3h_0+s)\ - \ s}{2}\ \ \ =\ \ \
\frac{3}{2}h_0 n\ +\ \frac{s(n-1)}{2}\ \ \ \ge\ \ \ \frac{3}{2}h_0 n\leqno{(6)}$$ This yields $n\le
2(h_0-1)/(\frac{3}{2}h_0-2)\le2$ since $h_0\ge2$. Plugging $n=2$ into (6) gives $h_0\le2$. Thus $h_0=2$ and
the inequalities in (6) are equalities, which implies $s=0$. This is case (I) in the Lemma.

\medskip\noindent
{\bf Case 4:} \ $h_0=1:$ \;  If $s=0$ then $g=1$, contrary to assumption. Thus $s>0$, and so (Hyp) gives
$r-3=s$. Thus $r\ge4$. Plugging $s=r-3$ into (5) gives
$$2n \ \ \ \ge\ \ \ \frac{r(n-1)\ +\ 3}{2}\leqno{(7)}$$
which implies $r\le\ 4+\frac{1}{n-1}$. If $n=2$ then we are in case (II) resp. (III).a of the Lemma. Assume
now $n>2$. Then $r=4$, $s=1$. If all $c_i=2$ then (7) would have to hold with equality, which is not
possible. Thus $c_4\ge3$ (by (4)).

Let $p_{i_0}$ be the image in $X/G$ of $q_1$. If $i_0<4$ then Riemann-Hurwitz yields $2n\ \ge \ \
\frac{3n-1}{2}+\frac{2}{3}n$, hence $n=3$ and $c_1=c_2=c_3=2$, a contradiction (since $n$ is divisible by
those $c_i$ with $i\ne i_0$). Thus $i_0=4$. If $n=3$ then $c_1=c_2=c_3=3$ contradicting Riemann-Hurwitz. Thus
$n>3$. Riemann-Hurwitz yields $2n\ \ge \ 3\ \frac{n}{2}+ \frac{2}{3}(n-1)$, hence $n=4$, $c_1=c_2=c_3=2$ and
$c_4=3$. Thus $\g_1$, $\g_2$,$\g_3$ act as double transpositions on the 4-set $G/H$, and $\g_4$ acts as
3-cycle. However, three double transpositions in $S_4$ cannot have a 3-cycle as their product (since the
double transpositions together with the identity form a normal subgroup). Contradiction.

\medskip\noindent
{\bf Case 5:} \ $h_0=0:$ \;  If $r\le2$ or $s\le2$ then $g=0$, contradiction. Thus we have $r\ge3$ and
$s\ge3$. Hence (Hyp) gives $r=s$. Plugging this into (5) gives $r\le4$.

\medskip\noindent
{\bf Case 5.a:} \ $r=4=s:$ \;  Here (5) holds with equality, hence the $q_j$'s are unramified over $X/G$,
and all other points of $X/H$ over some $p_i$ have ramification index $2=c_i$. If all $c_i=2$ then (1) gives
$g=1$, contradiction. Thus some $p_i$, say $p_{i_0}$, is unramified in $X/H$. If two $p_i$'s are unramified
in $X/H$ then $s\ge2n$, hence $n\le s/2=2$; this is case III.b in the Lemma.

Now suppose $p_i$ is unramified in $X/H$ only for $i=i_0$. Then $c_i=2$ for $i\ne i_0$, and $n$ is even
(since not each $p_i$, $i\ne i_0$, can have a $q_j$ over it); also, $n\le s=4$ and $n\ne2$ (because $X/H\to
X/G$ has 3 branch points). Hence $n=4$ and the $\g_i$ with $i\ne i_0$ act as double transpositions on the
4-set $G/H$ (and $\g_{i_0}$ acts trivially). This is case III.c.

\medskip\noindent
{\bf Case 5.b:} \ $r=3=s:$ \;  This is the last case to be considered. Now $X/H$ and $X/G$ have genus
zero, and the cover $X/H\to X/G$ (resp., $X\to X/H$) is ramified at 3 points of $X/G$ (resp., of $X/H$). This
gives the following condition on the action of $\g_1,\g_2,\g_3$ on the $n$-set $G/H$: Apart from exactly
three orbits of $\g_1,\g_2, \g_3$, all orbit lengths of $\g_i$ equal $c_i$; and if one of the three
exceptional orbits belongs to $\g_i$ then its length is $<c_i$. From (1) we also get
$\frac{1}{c_1}+\frac{1}{c_2}+\frac{1}{c_3}<1$. Finally, we have the genus zero condition
$2(n-1)=\sum_{i=1}^3\ \mbox{Ind}(\g_i)$. Triples of permutations with these properties (and with product 1,
and generating a transitive group) have been classified by Malle \cite{M}: His Table 2 gives the cycle types
of the three permutations plus the group they generate (plus the $c_i$). This is part (IV) of the Lemma.
(Actually, Malle's situation is slightly more general than ours, since he doesn't require $c_i$ to be the
order of $\g_i$; thus two of his cases, where $n=10$ and $n=18$, don't occur here). \endofproof

\section{The full automorphism group} \label{fullsection}
Now we change perspective and fix the finite group $H$. An $H$-curve $X$ is called {\bf generic} if the
following holds: The curve $X/H$, together with its points $q_1,...,q_s$ that are ramified in $X$, defines a
generic point of $\M_{h,s}$. Here $h$ is the genus of $X/H$, and \lq generic\rq \ means generic over
$\overline\bQ$ (the algebraic closure of $\bQ$).

Let $\Aut(X)$ denote the automorphism group of the curve $X$ (without regard of $H$-action).

\begin{Theorem} \label{thm1}
Suppose $X$ is a generic $H$-curve of genus $\ge 2$. Then $[\Aut(X): H]\le 24$. 
Let $h$ be the genus of $X/H$, and $s$ the number of points of $X/H$ ramified in $X$. Then $\Aut(X)=H$ unless $(h,s)$ equals $(2,0)$, $(1,2)$, $(1,1)$, $(0,4)$ or $(0,3)$.

Then there is an integer $n\le24$ such that $[\Aut(Y): H]\le n$ for all curves $Y$ corresponding to a dense open subset  of the locus $\M(g,H,\D)$. 
This is true with $n=1$ unless one of the following holds: If $h>0$ we have $[\Aut(X): H]\le 2$, and in the case $(0,4)$ we have $H$ normal in
$\Aut(X)$ and $[\Aut(X): H]\le 4$.
\end{Theorem}

\noindent \Proof Let $G=\Aut(X)$. Let $X$ be of type $(g,H,\D)$ when viewed as $H$-curve, and of type
$(g,G,\C)$ when viewed as $G$-curve.

Note that all maps and spaces in section \ref{section2} are defined over $\overline\bQ$. Hence the condition
that $X$ is a generic $H$-curve means that its corresponding point in $\M_g$ is a generic point of a
component of the locus $\M(g,H,\D)$. This point lies in the sublocus $\M(g,G,\C)$, hence the dimensions are
the same: \ $\d(g,G,\C)=\d(g,H,\D)$. Now the claim follows from Lemma \ref{Lemma2}. \endofproof

\begin{Lemma} \label{Lemma3}
\ (i) Let $\G$ be a group with a presentation on generators $\s_1,...\s_6$ and relations $$\s_1^2\ =...=\
\s_6^2\ =\ \ 1\ \ =\ \s_1\cdots\s_6$$ Let $\Pi$ be the kernel of the map $\G\to\{\pm1\}$, $\s_i\mapsto-1$.
Then $\Pi$ is generated by
$$\begin{array}{cc}
\a_1\ =\ \s_2\s_1\ \ \ \ \ \ & \a_2\ =\ \s_5\s_4 \\
\b_1\ =\ \s_2\s_3\ \ \ \ \ \ & \b_2\ =\ \s_5\s_6
\end{array}
$$
which gives a presentation for $\Pi$ with the single relation
$$ [\a_1,\b_1] \ [\a_2,\b_2] \ \ = \ \ 1$$
We have $\G=\ \Pi\xs\langle \s_2\rangle$ where $\s_2$ acts on $\Pi$ via
$$\begin{array}{cc}
\a_1 \mapsto \ \a_1^{-1}\ \ \ \ \ \ &
\a_2\mapsto \ \ (\a_2^{-1})^{\a_2\b_2\a_1^{-1}\b_1^{-1}} \\
\b_1\mapsto \ \b_1^{-1}\ \ \ \ \ \ & \b_2\mapsto \ \ (\b_2^{-1})^{\a_2\b_2\a_1^{-1}\b_1^{-1}}
\end{array} \leqno{(9)}$$
(ii)\ Let $\G=\ \langle \s_1,...\s_5\rangle$ subject to the relations
$$\s_1^2\ =...=\ \s_4^2\ =\ \ 1\ \ =\ \s_1\cdots\s_5$$ Let $\Pi$ be the
kernel of the map $\G\to\{\pm1\}$, $\s_i\mapsto-1$ for $i=1,...,4$ and $\s_5\mapsto1$. Then $\Pi$ is
generated by
$$\begin{array}{cc}
\a\ =\ \s_2\s_1,\ \ \ \ \ & \g_1\ =\ \s_5^{\s_4} \\
\b\ =\ \s_2\s_3,\ \ \ \ \ & \g_2\ =\ \s_5
\end{array}
$$
which gives a presentation for $\Pi$ with the single relation
$$ [\a,\b] \ \g_1\g_2 \ \ = \ \ 1$$
Now $\s_2$ acts on $\Pi$ via
$$\begin{array}{cc}
\a \mapsto \ \a^{-1}\ \ \ \ \ \ &
\g_1\mapsto \ \ \g_2^{\a^{-1}\b^{-1}} \\
\b\mapsto \ \b^{-1}\ \ \ \ \ \ & \g_2\mapsto \ \ \g_1^{\b^{-1}\a^{-1}}
\end{array} \leqno{(10)}$$

(iii)\ Let $\G=\ \langle \s_1,...\s_4\rangle$ subject to the relations
$$\s_1^2\ =\ \s_2^2\ =\ \s_3^2\ =\ \ 1\ \ =\ \s_1\cdots\s_4$$
Let $\Pi$ be the kernel of the map $\G\to\{\pm1\}$, $\s_i\mapsto-1$ (resp., $\s_3\mapsto(-1,-1)$, the normal
subgroup of $\G$ generated by the conjugates of $\s_4$). Then $\Pi$ is generated by
$$ \begin{array}{ccccc}
\begin{array}{cc}
\a\ =\ \s_2\s_1\ \ \ \ \ & \g\ =\ \s_4^2 \\
\b\ =\ \s_2\s_3\ \ \ \ \ &
\end{array} & \ \ \ \ & \mbox{resp.,} &\ \ \ \ &
\begin{array}{cc}
\g_1\ =\ \s_4\ \ \ \ \ & \g_3\ =\ \s_4^{\s_2} \\
\g_2\ =\ \s_4^{\s_1},\ \ \ \ \ & \g_4\ =\ \s_4^{\s_3}
\end{array}
\end{array}
$$
which gives a presentation for $\Pi$ with the single relation
$$ [\a,\b] \ \g \ = \ 1 \ \ \ \ \ \mbox{resp.,} \ \ \ \ \
\g_1\cdots \g_4\ = \ 1$$ Now $\s_2$ acts on $\Pi$ by inverting $\a$ and $\b$ (resp., $\s_2$ switches
$\g_1,\g_3$ and switches $\g_2,\g_4$, and $\s_1$ switches $\g_1,\g_2$ and maps $\g_3\mapsto\g_4^{\g_1}$,
$\g_4\mapsto\g_3^{\g_2^{-1}}$, and $\s_3$ switches $\g_1,\g_4$ and maps $\g_2\mapsto\g_3^{\g_4}$,
$\g_3\mapsto\g_2^{\g_1^{-1}}$).
\end{Lemma}

\noindent \Proof (i) All straightforward computation. One computes that the relation among the $\a_j$, $\b_j$
implies the relations among the $\s_i$, hence the former yields a presentation for $\Pi$ since $\G=\
\Pi\xs\langle \s_2\rangle$. The proof of (ii) and (iii) is similar. 
\endofproof

We say that $\a_1,\b_1,...,\a_{h},\b_{h},\g_1,...,\g_s$ form a {\bf genus $g$ generating system } of $H$ if
these elements generate $H$ and satisfy the basic relation (2), where $g$ is given by Riemann-Hurwitz
\[
\frac{2\ (g-1)}{|H|}\ \ \ = \ \ \ 2\ (h-1)\ +\ \sum_{i=1}^s\ \left( 1-\frac{1}{\mbox{ord} (\g_i)} \right)
\leqno{(8)}
\]
 Call a genus $g$ curve $X$ {\bf exceptional} if $X/\Aut(X)$ has genus $0$ and exactly three of
its points are ramified in $X$.

\begin{Theorem} \label{thm2} Fix $g\ge2$. A finite group $H$ is isomorphic
to the full automorphism group of a non-exceptional genus $g$ curve if and only if it has a genus $g$
generating system $\a_1,\b_1,...,\a_{h},\b_{h},\g_1,...,\g_s$ with $(h,s)\ne (0,3)$ satisfying:
\begin{enumerate}
\item[(a)]\ \ If $(h,s)=(2,0)$ (resp., $(1,2)$ ) then there is no automorphism of $H$ acting on the
generators via (9) (resp., (10) ). \item[(b)]\ \ If $(h,s)=(1,1)$ (resp., $(0,4)$ ) then there is no
automorphism of $H$ inverting $\a_1$ and $\b_1$ (resp., acting on the generators like $\s_1$, $\s_2$ or
$\s_3$ in (iii) of the Lemma).
\end{enumerate}
\end{Theorem}

For the proof we need the following set-up: Let $Y$ be a curve of genus $h$, and let $q_1,...,q_s$ be
distinct points of $Y$ such that the pair $(h,s)$ is as in (a) or (b). Then there is a Galois cover $Y\to\P$
such that the group $A=\Aut(Y/\P)$ permutes $q_1,...,q_s$; further, $|A|=2$ unless $(h,s)=(0,4)$ in which
case $A$ is a Klein 4 group whose non-trivial elements act as double transpositions on $\{q_1,...,q_4\}$. If
$(Y,q_1,...,q_s)$ corresponds to a generic point of $\M_{h,s}$ then $A$ is the full stabilizer of
$\{q_1,...,q_s\}$ in $\Aut(Y)$.

Let $\{p_1,...,p_r\}$ be the subset of $\P$ consisting of the points ramified in $Y$ plus the images of
$q_1,...,q_s$. Let $t$ be the identity function on $\P$, so that $\bC(t)$ is the function field of $\P$. View
$\bC(t)$ as a subfield of $\bC(Y)$ naturally. Let $\tilde\M$ be a maximal (algebraic) extension of $\bC(t)$
unramified outside the places $t=p_1,...,p_r$, and containing $\bC(Y)$. Its Galois group $\tilde\G$ is a
profinite group generated by elements $\tilde\s_1,...,\tilde\s_r$ subject to the single relation
$\tilde\s_1\cdots\tilde\s_r=1$; here $\tilde\s_i$ can be chosen to be a distinguished inertia group generator
over the place $t=p_i$. Let $\Omega$ be the smallest (closed) normal subgroup of $\tilde\G$ containing those
$\tilde\s_i^2$ for which there is no $q_j$ over $p_i$. Then $\G:=\tilde\G/\Omega$ is the Galois group of
$\M:=\tilde\M^\Omega$ over $\bC(t)$, and $\bC(t)\subset \bC(Y)\subset \M$. Let $\Pi=G(\M/\bC(Y))$, and let
$\s_1,...,\s_r$ be the images of the $\tilde\s_i$ in $\G$. Then $\G$ (resp., $\Pi$) is the profinite
completion of the group denoted by the same symbol in the Lemma; the cases $(h,s)=\ (2,0)$, $(1,2)$, $(1,1)$
and $(0,4)$ correspond to cases (i), (ii) and the two cases in (iii), respectively.

\medskip\noindent
\Proof (of the Theorem) \ First we prove that the conditon is necessary. So assume $H=\Aut(X)$ for a
non-exceptional genus $g$ curve $X$. Let $h$ be the genus of $Y:=X/H$, and $q_1,...,q_s$ the points of $Y$
ramified in $X$. It is well-known that then $H$ has a genus $g$ generating system with parameters $h,s$. So
we may assume that $(h,s)$ is as in (a) or (b). Then we are in the set-up discussed before this proof, with
$\bC(t)\subset \bC(Y)\subset \bC(X)\subset \M$. The generators of $\Pi$ from the Lemma yield generators of
$H$ via the natural surjection $\Pi=G(\M/\bC(Y))\to G(\bC(X)/\bC(Y))=H$. It remains to show that these
generators satisfy the condition in (a) and (b).

Let $f$ be the automorphism of $\Pi$ induced by conjugation action of $\s_2$ unless $(h,s)=(0,4)$ when $f$
could be induced by any of $\s_1,\s_2,\s_3$. If $H$ has an automorphism $\bar f$ making the following diagram
commutative
\[
\begin{matrix}\Pi & \mathop{\longrightarrow}\limits^{f} & \Pi\cr
\Big\downarrow & & \Big\downarrow\cr H & \mathop{\longrightarrow}\limits^{\bar f} & H\cr
\end{matrix}
\]
 then the kernel
$\Gamma$ of $\Pi\to H$ is $f$-invariant. Then the normalizer $\Delta$ of $\Gamma$ in $\G$ is strictly bigger
than $\Pi$; thus $\Delta/\Gamma$ is strictly bigger than $\Pi/\Gamma\cong H$. But $\Delta/\Gamma$ acts
faithfully on $X$ (since $\Gamma=G(\M/\bC(X))$. This contradicts the assumption $H=\Aut(X)$. Thus there is no
such commutative diagram, which implies (a) and (b).

Assume now $H$ has a genus $g$ generating system as in the Theorem. It is well-known that then there exists a
generic $H$-curve $X$ of genus $g$ with parameters $h$, $s$. If $(h,s)$ is not as in (a) or (b) then
$H=\Aut(X)$ by Theorem \ref{thm1} and we are done. Now assume $(h,s)$ is as in (a) or (b). Choose data
$(Y,q_1,...,q_s)$ as in the discussion before this proof, corresponding to a generic point of $\M_{h,s}$.
Mapping the generators of $\Pi$ from the Lemma to the given generators of $H$ defines a surjection $\Pi\to
H$. Let $K$ be the fixed field (in $\M$) of the kernel of this surjection. Then $K=\bC(X)$ for a generic
$H$-curve $X$ of genus $g$. Thus $H$ is normal in $G:=\Aut(X)$ by Theorem \ref{thm1}. Hence $G/H$ embeds into
$A$, and so $\bC(t)\subset K^G\subset \bC(Y)$. Thus the map $\Pi\to H$ extends to a map from a group between
$\Pi$ and $\G$ onto $G$. Thus if $G\ne H$ then there is an element of $G$ that acts on the given generators
of $H$ in the same way that $\s_2$ (resp., $\s_1$, $\s_2$ or $\s_3$ in case $(h,s)=(0,4)$ ) acts on the
generators of $\Pi$ (see Lemma). But this is excluded by the condition in (a) resp. (b). Hence $G=H$ and we
are done.\endofproof

\begin{Remark} The proof shows more: Given $g$, $H$ and a signature $\e$
(resp., a ramification type $(g,H,\E)$ ), the existence of a generating system of $H$ as in the Theorem with
$\mbox{ord}(\g_i)=e_i$ (resp., $\g_i\in E_i$) is equivalent to the existence of a non-exceptional genus $g$
curve $X$ with $\Aut(X)\cong H$ such that the resulting $H$-curve is of given signature (resp., type).
\end{Remark}

\part{Classification of automorphism groups}

Recall that a group $G$ acts faithfully on a genus $g$ curve if and only if it has a genus $g$ generating
system (see (8) above). For $g$ up to 48, all such groups and the signatures of all their genus $g$
generating systems have been listed by Breuer \cite{Br}. More precisely, for each $g\le48$, he produced a
list containing all {\bf signature-group pairs} in genus $g$, i.e., pairs consisting of a group $G$ together
with the signature of a genus $g$ generating system of $G$.

If $G$ acts on $X_g$ then so does each subgroup of $G$. This shows that Breuer's lists have to be long, and
contain some redundancy. Part II of our paper was written in order to eliminate those signature-group pairs
that do not yield the full automorphism group of a curve. This can also be done by the results of \cite{Si}
and \cite{Ri} (as we learned later). It turns out that the larger $g$ gets, the larger is the ratio of
entries in Breuer's lists that do occur as full automorphism group in genus $g$. That can already be seen
from the fact (see Lemma \ref{Lemma2}) that if a signature-group pair does not yield the full automorphism
group of a curve, then its $\delta$-invariant (dimension of corresponding locus in $\M_g$) is at most 3.

For small $g$, a relatively large portion of those groups do not occur as full automorphism group in genus
$g$. Among those that do occur, we distinguish those that occur for a particularly simple class of curves:
Call a group homocyclic if it is a direct product of isomorphic cyclic groups. A cyclic (resp., homocyclic)
cover $X\to\P$ is defined to be a Galois cover with cyclic (resp., homocyclic) Galois group $C$. We call it a
{\bf normal cyclic (resp., homocyclic) cover of} $\P$ if $C$ is normal in $G=\Aut(X)$. Then $\bar G:=G/C$
embeds as a finite subgroup of $PSL(2,\bC)$ and it is easy to write down an equation for $X$ (e.g., in the
cyclic case, $X$ has an equation of the form $y^c=f(x)$, where $c=|C|$). To illustrate this, we work out the
equations in the genus 3 case. The hyperelliptic case (cyclic with $|C|=2$) is the simplest case and can be
done as in the genus 2 case (see \cite{ShV}); see also \cite{Bra} for the cyclic case (in arbitrary
characteristic).

\section{Genus 3}
For $g=3$, there are 49 signature-group pairs. Applying Theorem \ref{thm2} and Lemma \ref{Lemma2} we obtain
that 23 of them correspond to $G$-curves having $G$ as full automorphism group. For each of these 23
signature-group pairs, we find that the corresponding curves form an irreducible locus in $\M_3$, and we
compute the equation of a general curve in that locus. There is only two groups (namely $C_2$ and $V_4$) that
occur with more than one signature, and as soon as $|G|>2$, then the orbit genus $g_0$ (genus of $X_3/G$) is
always $0$.

Each row in Table 1,2 or 3 corresponds to an irreducible family of curves $X_3$ with automorphism group $G$
(given in the first column). We let $K$ denote the function field of $X_3$ over the algebraically closed
field $k$ of characteristic $0$ (and identify $G$ with $\Aut(K/k)$). We display an equation for $K$. For
almost all values (Zariski-open set) of the parameters in these equations, the equation describes a genus 3
curve with $G$ as full automorphism group. For the exceptional values, the automorphism group may be larger,
or the curve may degenerate. (We do not specify these exceptional values). Most of these equations agree with
those found in \cite{KuKo1}, \cite{KuKo2}, but they missed two cases: The hyperelliptic case with $\bar
G=D_{12}$, and the 3-dimensional hyperelliptic case with $\bar G=V_4$. Further, they display an equation for
the hyperelliptic case with $\bar G=S_4$ which actually belongs to the case $\bar G=D_{12}$. Some of the loci
they describe with more parameters than the dimension requires. The cases with $g_0=1$ and group of order 3
or 4 should be omitted in their list, because the full automorphism group of these curves is larger (as shown
by our Theorem \ref{thm2} above). They note this for two of these cases in the Errata at the end of their
paper, but again missed one case.

Now we show how to obtain these equations quickly from the group-theoretic data. Starting point is Breuer's
data in genus 3 (see above). We used Part II to eliminate those signature-group pairs that do not yield the
full automorphism group of a genus 3 curve. Now we discuss the remaining signature-group pairs.

\subsection{Homocyclic non-hyperelliptic case}
We show how to derive the information in Table 2.

\smallskip\noindent
$G=C_4^2\xs S_3$: \ The fixed field of $C_4^2$ is of genus $0$, call it $k(u)$. The extension $K/k(u)$ is
ramified at 3 places of $k(u)$, say $u=0,1,\infty$, and the corresponding inertia groups are cyclic of order
4. The fixed fields of these 3 inertia groups extend $k(u)$ with 2 ramified places each, hence they are of
genus $0$; call them $k(x)$, $k(y)$ and $k(z)$. We may assume $x^4=u-1$, $y^4=u$. Thus $K=k(x,y)$ with
$y^4=u=x^4+1$. The equation $y^4=x(x^2-1)$ given by \cite{KuKo1} is also correct. (Indeed, the four roots of
$x^4+1$ have $j$-invariant 1728, same as $\infty$ plus the 3 roots of $x(x^2-1)$).

\medskip\noindent
$G=(48,33)$: Here $G$ is a central extension of $C_4$ by quotient $A_4$. Let $H\cong C_3$ be a Sylow
3-subgroup of $G$. The fixed field of $H$ is a $k(x)$, and $K/k(x)$ is ramified at 5 places of $k(x)$. Let
$h_1,...,h_5$ be the corresponding distinguished inertia group generators in $H$. Their product is 1, hence
we may assume $h_1=...=h_4=h_5^{-1}$. Take the place with generator $h_5$ to be $x=\infty$, and the other 4
ramified places to be the roots of $f(x)$. Then $K=k(x,y)$ with $y^3=f(x)$. Further, the $C_4$ acts on $k(x)$
permuting the 4 roots of $f(x)$. We can choose $x$ such that $C_4$ acts on $k(x)$ by the maps $x\mapsto \zeta
x$, $\zeta^4=1$. Then $f(x)=x^4-b$. We can get each $b\ne0$ by further normalization of $x$. This agrees with
\cite{KuKo1}.

\medskip\noindent
$G=(16,13)$: Here $G$ is a central extension of $H\cong C_4$ by quotient $V_4$. The fixed field of $H$ is a
$k(x)$, and $K/k(x)$ is ramified at 4 places $p_1,...,p_4$ of $k(x)$. Let $h_1,...,h_4$ be the corresponding
inertia elements in $H$. The Klein 4 group $G/H$ permutes $p_1,...,p_4$. If this action were not transitive,
then $p_1,...,p_4$ could not be 4 general points in $\P$, contradicting the fact that the covers under
consideration form a 1-dimensional family. Thus $G/H$ permutes $p_1,...,p_4$ transitively, hence permutes
$h_1,...,h_4$ transitively (by conjugation). Thus $h_1=...=h_4$ since $H$ is central in $G$. (Alternatively,
this can be shown by the algorithm before Lemma \ref{Lemma2}). Hence $K=k(x,y)$ with $f$ of degree 4. The
generic curve of this type has automorphism group of order 16, because there is always a Klein 4 group in
Aut$(k(x))$ permuting the 4 roots of $f$.

\medskip\noindent
$G=C_9,C_6$ or $C_3$: Let $H$ be the subgroup of $G$ of order 3. The fixed field of $H$ is a $k(x)$, and
$K/k(x)$ is ramified at 5 places $p_1,...,p_5$ of $k(x)$. Let $h_1,...,h_5$ be the corresponding inertia
elements. They have product 1, hence we may assume $h_1=...=h_4=h_5^{-1}$. The group $G/H$ permutes
$p_1,...,p_5$, fixing $p_5$. Take $p_1,...,p_5$ to be the places $x=0,1, s,t,\infty$, respectively. Then
$K=k(x,y)$ with $y^3=x(x-1)(x-s)(x-t)$. Let $\alpha$ be the automorphism of $k(x)$ induced by a generator of
$G/H$. If it has order 3 then we may assume that $x=0$ is the other fixed point of $\alpha$ (besides
$x=\infty$); then $\alpha(x)=\zeta_3x$ and the equation becomes $y^3=x(x^3-1)$. Now assume $\alpha$ has order
2. We may assume it switches $0$ and $1$. Then $\alpha(t)=1-t$, again as claimed.

\medskip\noindent
$G=S_4, D_8$ or $V_4$: Let $H$ a normal Klein 4-group in $G$. The fixed field of $H$ is a $k(w)$, and
$K/k(w)$ is ramified at 6 places $p_1,...,p_6$ of $k(w)$. Let $h_1,...,h_6$ be the corresponding inertia
elements. If they don't comprise all involutions in $H$ then there is $h\in H$ that occurs 4 times among the
$h_i$; then $K^h$ extends $k(w)$ with 2 branch points, hence $K^h$ has genus $0$ and so $K$ is hyperelliptic.
We excluded that case. Thus we may assume $h_1=h_2$, $h_3=h_4$, $h_5=h_6$ and $\{h_1,h_3,h_5\}=
H\setminus\{1\}$.

Let $u,v$ be independent transcendentals over $k$ with $w=u/v$. The places of $k(w)$ correspond to
homogeneous degree 1 polynomials in $u,v$, modulo scalar multiples. Let $P,Q,R$ be homogeneous degree 2
polynomials in $k[u,v]$ such that the degree 1 factors of $P$ (resp., $Q$, resp., $R$) correspond to
$p_1,p_2$ (resp., $p_3,p_4$, resp., $p_5, p_6$). The quotients $P/Q$, $R/Q$, $P/R$ are naturally elements of
$k(w)$. Then $K^{h_3}=k(w,x)$ with $x^2=P/Q$, and $K^{h_1}= k(w,y)$ with $y^2=R/Q$, and $K^{h_2}=k(w,z)$ with
$z^2=P/R$. Thus $K=k(w,x,y)$. We claim $K =k(x,y)$. Assume wrong. Then $[K:k(x,y)]=2$ since $w$ satisfies the
equation $P(w,1)-x^2Q(w,1)=0$ over $k(x,y)$. Let $g_0$ be the generator of $G(K/k(x,y))$. Then $g_0$ is not
in $H$, since none of the $h_i$ fixes both $x$ and $y$. But each $h_i$ maps $x$ (resp., $y$) to $\pm x$
(resp., $\pm y$), hence $H$ leaves $k(x,y)$ invariant and thus centralizes $g_0$. Hence $<H,g_0>\cong C_2^3$,
contradiction (since $G$ contains no such subgroup). This proves $K=k(x,y)$.

The 6 homogeneous polynomials $P^2, Q^2, R^2, PQ, PR, QR$ of degree 4 are linearly dependent in $k[u,v]$.
This gives a non-trivial relation $$aPR\ +\ bPQ\ +\ cQR\ +\ d P^2\ +\ eQ^2\ +\ fR^2\ = \ 0$$ If $f=0$ then
dividing by $Q^2$ gives an equation for $y$ over $k(x)$ of $y$-degree $\le2$ and so $K$ is hyperelliptic or
of genus $0$ -- contradiction. Thus $f$ (and similarly $d,e$) are non-zero. Replacing $P$, $Q$, $R$ by
suitable scalar multiples we get a relation $$ P^2\ +\ Q^2\ +\ R^2\ +\ aPR\ +\ b PQ\ +\ c QR \ = \ 0$$ with
$a,b,c\in k$. Dividing by $Q^2$ gives $$x^4\ +\ 1\ +\ y^4\ +\ ax^2y^2\ +\ bx^2\ +\ cy^2 \ = \ 0$$ This
settles the case $G=V_4$.

Suppose now $G=D_8$. Then $G/H$ induces an involutory automorphism $\alpha$ of $k(w)$ switching (without
loss) the sets $\{p_1,p_2\}$ and $\{p_5,p_6\}$. We may assume $\alpha(w)=1/w$. Extend it to $k(u,v)$ by
$\alpha(u)=v$, $\alpha(v)=u$. Then $\alpha$ fixes $Q$ and switches $P,R$ up to scalar multiples. Since
$\alpha^2=\mbox{id}$, we have $\alpha(Q)=\pm Q$. Applying $\alpha$ to the relation between $P,Q,R$ yields
that $\alpha(P)=\pm R$. Replacing $R$ by $-R$ if necessary we may assume $\alpha(P)=R$ (resp.,
$\alpha(P)=-R$) if $\alpha(Q)=Q$ (resp., $\alpha(Q)=-Q$). Then $b=c$. Finally, suppose $G=S_4$. Then there is
another automorphism like $\alpha$, fixing $R$ and switching $P,Q$ up to scalar multiples. Replacing $Q$ by
$-Q$ if necessary yields $a=b=c$.

\subsection{Hyperelliptic case}
Here we derive the information in Table 3.

In the hyperelliptic case, $G$ has a central subgroup $C$ of order 2 whose fixed field in $K$ is a $k(x)$.
There are 8 places of $k(x)$ that ramify in $K$, say $x=a_1,...,a_8$. Then $\bar G=G/C$ embeds into
$PGL(2,k)$ as a subgroup permuting $a_1,...,a_8\in\P$. Let $\sigma_m$ (resp., $\tau$) be the element of
$PGL(2,k)$ mapping $x$ to $\zeta_m x$ (resp., $1/x$), where $\zeta_m$ is a primitive $m$-th root of 1 in $k$.
Then each subgroup of $PGL(2,k)$ isomorphic to $C_m$ resp., $D_{2m}$ is conjugate to $<\sigma_m>$ resp.,
$<\sigma_m, \tau>$. We identify $C_m$ (resp., $D_m$) with $<\sigma_m>$ (resp., $<\sigma_m,\tau>$). In
particular, $D_2=V_4=<\sigma_2,\tau>$.

\smallskip\noindent
$\bar G=S_4$: Let $i\in k$ with $i^2=-1$. We may assume $\bar G$ contains $D_2$ and
$\rho:x\mapsto(x+i)/(x-i)$. ($\rho$ is of order 3, permuting the 6 fixed points of the involutions in $V_4$).
The two fixed points of $\rho$ are the roots of $X^2-(i+1)X-i=0$. The group $<V_4,\rho>\cong A_4$ has exactly
one orbit of length 8 on $\P$, containing those fixed points of $\rho$. Hence the orbit consists of the roots
of $(X^2-(i+1)X-i) (X^2+(i+1)X-i)(-iX^2-(i+1)X+1)(-iX^2+(i+1)X+1) \ = \ X^8+14X^4+1$.

\medskip\noindent
$\bar G=D_{16}, D_{12},C_7$: In the dihedral case, $\bar G$ leaves exactly one set of length $8$ in $\P$
invariant, consisting of the roots of $X^8-1$, resp., the roots of $X^6-1$ plus $0$ and $\infty$. In the
third case, we may assume $a_1=1$ by a coordinate change; then $a_1,...,a_8$ comprise the roots of $X^7-1$
plus either $0$ or $\infty$. (The latter two choices are conjugate under the map $\tau$). The case $\bar
G=D_{12}$ is missing in \cite{KuKo1} and the corresponding equation $y^2=x(x^6-1)$ is erroneously associated
with the case $\bar G=S_4$.

\medskip\noindent
$\bar G=D_8,S_3$: Here $a_1,...,a_8$ comprise a regular orbit of $D_8$ resp., $0$ and $\infty$ plus a regular
orbit. Thus in the first case, $a_1,...,a_8$ are the roots of $(X^4-t)(X^4-1/t)= X^8-(t+1/t)X^4+1$. The other
case is similar.

\medskip\noindent
$\bar G=V_4$: Here $a_1,...,a_8$ either comprise two regular orbits of $\bar G$, or a regular orbit plus two
of three orbits of length 2. The orbits of length 2 are all conjugate under $S_4$, so we may take
$\{0,\infty\}$ and $\{\pm1\}$ as the two orbits of length 2. Now the claim follows as in the previous case.

\medskip\noindent
$\bar G=C_2$: Here $a_1,...,a_8$ either comprise 4 regular orbits, or $0$ and $\infty$ plus 3 regular orbits.
We may assume $a_1=1$ by a coordinate change.

\subsection{Remaining cases, collected in Table 1}

Here we take $X_3$ in its canonical embedding in the projective plane associated with the dual of the
3-dimensional (faithful) $G$-module $\Omega$ of holomorphic differential forms on $X_3$. Let $x,y,z$ be
coordinates of this module $\Omega^*$. Let $f(x,y,z)=0$ be the equation of $X_3$ in these coordinates. Then
$f$ spans a 1-dimensional fixed space for the representation of $G$ on the space of quartic homogeneous
polynomials in $x,y,z$.

To identify the representation of $G$ on $\Omega$, we use that the space of fixed points of $G$ on $\Omega$
has dimension equal to the genus $g_0$ of $X_3/G$ (see \cite{Br}); this fixed point space can be identified
with $\Omega(X_3/G)$).

\medskip\noindent
$G=L_3(2)$: This is the well-known Hurwitz group of order 168. The $(2,3,7)$ triple is weakly-rigid in $G$,
so there is exactly one associated curve found already by Klein.

\medskip\noindent
$G=S_3, g_0=0$: The group $S_3$ has only one faithful 3-dimensional representation with trivial fixed space:
It is the natural permutation representation tensored with the 1-dimensional sign-representation. Thus $G$
acts in the projective plane associated with $\Omega^*$ just by permuting the homogeneous coordinates $x,y,z$
(for suitable choice of such coordinates). Thus the quartic polynomial $f(x,y,z)$ defining $X_3$ is either
symmetric in $x,y,z$, or has the property that it is invariant under a cyclic permutation of $x,y,z$, but
changes sign under a transposition. The only possibility (up to scalar multiples) for the latter case is 
$
f=x^3y+y^3z+z^3x-y^3x-z^3y-x^3z,
$
 but this polynomial is not irreducible (product of linear factors). Hence $f$
is symmetric, thus is of the form 
\[
f=a(x^4+y^4+z^4)+b(x^2y^2+x^2z^2+y^2z^2)+c(x^2yz+y^2xz+z^2xy)+
d(x^3y+x^3z+y^3x+y^3z+z^3x+z^3y).
\]
 The centralizer of $S_3$ in $GL_3(\bC)$ consists of the matrices $tI+sJ$,
where $I$ (resp., $J$) is the identity matrix (resp., the all $1$'s matrix). By a coordinate change with such
a matrix we can make $d=0$, but not uniquely: Generically, there is always 4 triples $(a,b,c)$ (up to scalar
multiples) that give the same curve (up to isomorphism).

\medskip\noindent
$G=C_2, g_0=1$:\ We may assume that the involution in $G$ acts as $x\mapsto-x, y\mapsto y, z\mapsto z$. Thus
$f$ either contains only even powers of $x$ (i.e., $f$ is fixed by $G$) or $f$ contains only odd powers of
$x$ (i.e., $f$ is mapped to $-f$ by the involution in $G$). In the latter case, $f$ is reducible (product of
$x$ and a genus 1 equation). Thus $f$ contains only even powers of $x$, hence 
\[
f\ \ =\ \ a_0x^4\ +\ x^2\
(a_1y^2+a_2yz+a_3z^2)\ + \ a_4y^4+a_5y^3z+a_6y^2z^2+a_7yz^3+a_8z^4
\]
 By a linear change of the variables
$y,z$ we can make $a_1=1$ and $a_2=0$. (If $a_1=a_2=a_3=0$ then $X_3$ has an automorphism of order 4, namely
multiplying $x$ by $\zeta_4$, so we don't consider this case here). Further we can replace $x$ and $z$ by
scalar multiples to get the normalizations from Table 1. (If either $a_0=0$ or $a_7=a_8=0$ then the curve is
hyperelliptic).

\begin{table}[h] \label{genus 3.nothomo}
\caption{$\Aut(X_3)$ for $X_3$ not a normal homocyclic cover of $\P$}  

\begin{center}
\begin{tabular}{||c|c|c|c|c||}

\hline \hline &&&&\\ $G$ & signature & orbit& proj. equation & Group
\\ &&genus &&ID\\
\hline \hline &&&&\\ $L_3(2)$ & $(2,3,7)$ & 0
&$x^3y+y^3z+z^3x=0$ &(168,42) \\

\hline &&&&\\
$S_3$ & $(2,2,2,2,3)$ &0 &
$a(x^4+y^4+z^4)+b(x^2y^2+x^2z^2+y^2z^2)+$&(6,1)\\

&&&$c(x^2yz+y^2xz+z^2xy)=0$&\\

\hline &&&&\\ $C_2$ & $(2,2,2,2)$ &1&$ x^4+x^2(y^2+az^2) + by^4+cy^3z+dy^2z^2+eyz^3+gz^4=0$&(2,1)\\
&&& either $e=1$ or $g=1$ &\\ \hline \hline
\end{tabular}
\end{center}
\end{table}

\begin{table}[h] \label{genus 3.homo}
\caption{$\Aut(X_3)$ for normal homocyclic covers of $\P$, not hyperelliptic}  

\begin{center}
\begin{tabular}{||c|c|c|c|c|c||}
\hline \hline &&&&&\\ $G$ & $C$ &$ G/C$ & signature & equation & Group
\\ &&&&&ID\\ \hline


\hline &&&&&\\
$V_4$ & $V_4$ & $\{1\}$ & $(2^6)$ &$x^4+y^4+ax^2y^2+bx^2+cy^2+1=0$&(4,2) \\

$D_8$ & $V_4$ & $C_2$ & $(2^5)$ & take\ $b=c$& (8,3)\\

$S_4$ &$V_4$& $S_3$ & $(2,2,2,3)$ & take\ $a=b=c$ &(24,12) \\

$C_4^2\xs S_3$ &$V_4$ & $S_4$ & $(2,3,8)$ &
take\ $a=b=c=0$ \ or\ $y^4=x(x^2-1)$& (96,64)\\

\hline

$16$ & $C_4$ & $V_4$& $(2,2,2,4)$ &$y^4=x(x-1)(x-t)$&(16,13)
\\

$48$ &$C_4$& $A_4$& $(2,3,12)$ &$y^4=x^3-1$
&(48,33)\\

\hline

$C_3$ &$C_3$& $\{1\}$ &$(3^5)$ &$y^3=x(x-1)(x-s)(x-t)$&(3,1)\\

$C_6$ & $C_3$ & $C_2$&$(2,3,3,6)$ & take\ $s=1-t$& (6,2)\\

$C_9$ &$C_3$ &$C_3$ & $(3,9,9)$
&$y^3=x(x^3-1)$&(9,1)\\

\hline \hline

\end{tabular}
\end{center}
\end{table}

\begin{table}[] \label{genus 3.hyperelliptic}
\caption{$\Aut(X_3)$ for hyperelliptic $X_3$} 
\begin{center}
\begin{tabular}{||c|c|c|c|c|c||}
\hline \hline &&&&&\\ $G$ & $ G/C$ & signature & dim. of&equation & Group
\\ &&&locus&$y^2=$&ID\\ 
\hline \hline &&&&&\\
$48$ & $S_4$ &$(2,4,6)$ &0 & $x^8+14x^2+1$ &$(48,48)$\\
$32$ &$D_{16}$ &$(2,4,8)$ &0& $x^8-1$& $(32,9)$\\
$24$ &$D_{12}$ & $(2,4,12)$ &0&$x(x^6-1)$ & $(24,5)$\\
$C_{14}$ &$C_7$& $(2,7,14)$ &0&$x^7-1$ &$(14,2)$ \\
$16$ &$D_8$ &$(2,2,2,4)$&1&$x^8+ax^4+1$ &$(16,11)$ \\
$D_{12}$ &$S_3$& $(2,2,2,6)$ &1&$x(x^6+ax^3+1)$ &$(12,4)$\\
$C_2\times C_4$ &$V_4$ &$(2,2,4,4)$ & 1&$x(x^2-1)(x^4+ax^2+1)$& $(8,2)$\\
$C_2^3$ &$V_4$& $(2^5)$ &2&$(x^4+ax^2+1)(x^4+bx^2+1)$ &$(8,5)$\\
$C_4$ &$C_2$& $(2,2,2,4,4)$ &2&$x(x^2-1)(x^4+ax^2+b)$& $(4,1)$\\
$V_4$ &$C_2$& $(2^6)$ & 3&$(x^2-1)(x^6+ax^4+bx^2+c)$ &$(4,2)$\\
$C_2$ &$\{1\}$& $(2^8)$ & 5&$x(x-1)(x^5+ax^4+bx^3+cx^2+dx+e)$& $(2,1)$
\\ \hline\hline

\end{tabular}
\end{center}
\end{table}

\clearpage
\section{Curves of genus $\le10$ with large automorphism group}
In the genus 3 case, we were able to write out explicit equations for the curves with any given automorphism
group. This yields an explicit description of the corresponding loci in $\M_3$. For higher genus, we cannot
expect to obtain explicit equations. Still, computational group theory allows us to determine the dimensions
and number of components of these loci. This uses the general set-up of section 3.

The number of signature-group pairs grows quickly with the genus. E.g., in genus 10 there are already 174
signature-group pairs with $g_0 = 0$, and most of them yield the full automorphism group of a genus 10 curve.
So it would not be feasible to display all automorphism groups up to genus 10. Therefore, we only display the
{\bf large} groups Aut$(X_g)$ (in the sense of section 1.3), see Table 4. Surprisingly, their number remains
relatively small. They comprise the most interesting groups in each genus, and we avoid listing the many
group-signature pairs with group of order 2, 3 etc.

\subsection{The general set-up}
We return to the set-up of section 3. Let $\c=(c_1,...,c_r)$ be the signature of a genus $g$ generating
system of $G$. Let $\H(g,G,\c)$ be Hurwitz space parameterizing equivalence classes of $G$-curves of genus
$g$ and signature $\c$; i.e., $\H(g,G,\c)$ is the (disjoint) union of the spaces $\H(g,G,\C)$ with $\C$ of
signature $\c$ (see section 3). The map
$$\Phi:\ \ \H(g,G,\c) \ \ \to \ \ \M_{g}$$
forgetting the $G$-action is a finite morphism. Let $\M(g,G,\c)$ be its image (the locus of genus $g$ curves
admitting a $G$-action of signature $\c$). All components of $\M(g,G,\c)$ have dimension
$\delta=\delta(g,G,\c)$ (by Lemma 3.1). In particular, if $G$ is large then $g_0=0$ and $r=3,4$, hence
$\delta$ is $0$ or $1$.

Define $\H^*(g,G,\c)$ (resp., $\H^*(g,G,\C)$ with $\C$ a tuple of conjugacy classes, see section 3) to be the
union of all components $C$ of $\H(g,G,\c)$ (resp., $\H(g,G,\C)$ ) with the following property: There is at
least one point on $C$ such that the associated $G$-curve has $G$ as full automorphism group. Let
$\M^*(g,G,\c)$ (resp., $\M^*(g,G,\C)$ ) be the $\Phi$-image of the corresponding space $\H^*(...)$. Then
$\Phi$ is generically injective on $\M^*(g,G,\c)$, and so the spaces $\M^*(g,G,\c)$ and $\H^*(g,G,\c)$ have
the same number of components. In the case $g_0=0$, the spaces $\H(g,G,\c)$ coincide with the Hurwitz spaces
studied in \cite{FV}, \cite{V}, \cite{Buch}. Thus the components of $\H(g,G,\c)$ correspond to the braid
orbits of (genus $0$) generating systems of $G$ of signature $\c$, taken modulo $\Aut(G)$. These braid orbits
can be computed with the BRAID package \cite{MSV}.

\subsection{The table of large automorphism groups up to genus $10$}
In Table 4 we list all group-signature pairs $(G,\c)$ of genus $g$, where $4 \le g \le 10$, with the
following property: There exists a $G$-curve $X_g$ of genus $g$ and signature $\c$ such that $G$ is the full
automorphism group of $X_g$ and $G$ is large; i.e, \[ |G|> 4(g-1).\]
 More precisely, the rows of Table 4
correspond to the components of the loci $\M^*(g,G,\c)$ associated with these parameters. It turns out that
these loci are mostly irreducible, with only 6 exceptions listed in Table 5. In these exceptional cases, they
have two components, and these components are all of the form $M^*(g,G,\C)$, with $\C$ a tuple of conjugacy
classes (see 7.1). Thus the spaces $M^*(g,G,\C)$ are always irreducible in the situation of Table 4, and they
correspond bijectively to the rows of Table 4. In particular, duplicate rows occur iff the corresponding
locus $\M^*(g,G,\c)$ is reducible. The group $G$ is identified via its ID from the Small Group Library. In
the last column of Table 4, we also indicate the inclusion relations between components of dimension 0 and 1.
They can be computed by the algorithm given before Lemma \ref{Lemma2}.

\begin{table}[h] \label{genus 4-10}
      \caption{Components of the Hurwitz loci $\M^*(g,G,\C)$ for large $G$}
      \begin{center}
      \renewcommand{\arraystretch}{1.24}
      \begin{tabular}{||c|c|c|c||c|c|c|c||}
      \hline
      \hline
      $\#$ & Group ID & signature & contains & $\#$ & Group ID & signature & contains \\
      \hline
      \hline
      \multicolumn{8}{c}{Genus 4, $\delta = 0$} \\
      \hline
      \hline
      {\bf 1} & (120,34) & $(2, 4, 5)$ & &
      {\bf 2} & (72,42) & $(2, 3, 12)$ & \\
      {\bf 3}& (72,40) & $(2, 4, 6)$ & &
      {\bf 4} & (40,8) & $(2, 4, 10)$ & \\
      {\bf 5} & (36,12) & $(2, 6, 6)$ & &
      {\bf 6} & (32,19) & $(2, 4, 16)$ & \\
      {\bf 7} & (24,3) & $(3, 4, 6)$ & &
      {\bf 8} & (18,2) & $(2, 9, 18)$ & \\
      {\bf 9} & (15,1) & $(3, 5, 15)$ &&&&&\\
      \hline
      \hline
      \multicolumn{8}{c}{Genus 4, $\delta = 1$} \\
      \hline
      \hline
      {\bf 10} & (36,10) & $(2, 2, 2, 3)$ & 3 &
      {\bf 11} & (24,12) & $(2, 2, 2, 4)$ & 1, 2  \\
      {\bf 12} & (20,4) & $(2, 2, 2, 5)$ & 4 &
      {\bf 13} & (18,3) & $(2, 2, 3, 3)$ & 2, 5  \\
      {\bf 14} & (16,7) & $(2, 2, 2, 8)$ & 6 &&&& \\
      \hline
      \hline
      \multicolumn{8}{c}{Genus 5, $\delta = 0$} \\
      \hline
      \hline
      {\bf 1} & (192,181) & $(2, 3, 8)$ & &
      {\bf 2} & (160,234) & $(2, 4, 5)$ & \\
      {\bf 3} & (120,35) & $(2, 3, 10)$ & &
      {\bf 4} & (96,195) & $(2, 4, 6)$ & \\
      {\bf 5} & (64,32) & $(2, 4, 8)$ & &
      {\bf 6} & (48,14) & $(2, 4, 12)$ & \\
      {\bf 7} & (48,30) & $(3, 4, 4)$ & &
      {\bf 8} & (40,5) & $(2, 4, 20)$ & \\
      {\bf 9} & (30,2) & $(2, 6, 15)$ & &
      {\bf 10} & (22,2) & $(2, 11, 22)$ & \\
      \hline
      \hline
      \multicolumn{8}{c}{Genus 5, $\delta = 1$} \\
      \hline
      \hline
      {\bf 11} & (48,48) & $(2, 2, 2, 3)$ & 1, 4 &
      {\bf 12} & (32,43) & $(2, 2, 2, 4)$ & \\
      {\bf 13} & (32,28) & $(2, 2, 2, 4)$ & 1 &
      {\bf 14} & (32,27) & $(2, 2, 2, 4)$ & 2, 4, 5 \\
      {\bf 15} & (24,14) & $(2, 2, 2, 6)$ & 6 &
      {\bf 16} & (24,8) & $(2, 2, 2, 6)$ & 4 \\
      {\bf 17} & (24,13) & $(2, 2, 3, 3)$ & 3, 7 &
      {\bf 18} & (20,4) & $(2, 2, 2, 10)$ & 3, 8 \\

      \hline
      \hline
      \multicolumn{8}{c}{Genus 6, $\delta = 0$} \\
      \hline
      \hline
      {\bf 1} & (150,5) & $(2, 3, 10)$ & &
      {\bf 2} & (120,34) & $(2, 4, 6)$ & \\
      {\bf 3} & (72,15) & $(2, 4, 9)$ & &
      {\bf 4} & (56,7) & $(2, 4, 14)$ & \\
      {\bf 5} & (48,6) & $(2, 4, 24)$ & &
      {\bf 6} & (48,29) & $(2, 6, 8)$ & \\
      {\bf 7} & (48,15) & $(2, 6, 8)$ & &
      {\bf 8} & (39,1) & $(3, 3, 13)$ & \\
      {\bf 9} & (30,1) & $(2, 10, 15)$ & &
      {\bf 10} & (26,2) & $(2, 13, 26)$ & \\
      {\bf 11} & (21,2) & $(3, 7, 21)$ &&&&&\\
      \hline
      \hline
      \multicolumn{8}{c}{Genus 6, $\delta = 1$} \\
      \hline
      \hline
      {\bf 12} & (60,5) & $(2, 2, 2, 3)$ & 2 &
      {\bf 13} & (28,3) & $(2, 2, 2, 7)$ & 4 \\
      {\bf 14} & (24,12) & $(2, 2, 3, 4)$ & 2 &
      {\bf 15} & (24,8) & $(2, 2, 3, 4)$ & 3 \\
      {\bf 16} & (24,6) & $(2, 2, 2, 12)$ & 5 &
      {\bf 17} & (24,6) & $(2, 2, 3, 4)$ & 7 \\

      \hline
      \hline
      \end{tabular}
      \end{center}
      \end{table}

      \addtocounter{table}{-1}
      \begin{table}

      \caption{(Cont.)}
      \begin{center}
      \renewcommand{\arraystretch}{1.24}
      \begin{tabular}{||c|c|c|c||c|c|c|c||}
      \hline
      \hline
      $\#$ & Group ID & signature & contains & $\#$ & Group ID & signature &
contains \\

\hline
      \hline

      \multicolumn{8}{c}{Genus 7, $\delta = 0$} \\
      \hline
      \hline
      {\bf 1} & (504,156) & $(2, 3, 7)$ & &
      {\bf 2} & (144,127) & $(2, 3, 12)$ & \\
      {\bf 3} & (64,41) & $(2, 4, 16)$ & &
      {\bf 4} & (64,38) & $(2, 4, 16)$ & \\
      {\bf 5} & (56,4) & $(2, 4, 28)$ & &
      {\bf 6} & (54,6) & $(2, 6, 9)$ & \\
      {\bf 7} & (54,6) & $(2, 6, 9)$ & &
      {\bf 8} & (54,3) & $(2, 6, 9)$ & \\
      {\bf 9} & (48,32) & $(3, 4, 6)$ & &
      {\bf 10} & (42,4) & $(2, 6, 21)$ & \\
      {\bf 11} & (32,11) & $(4, 4, 8)$ & &
      {\bf 12} & (32,10) & $(4, 4, 8)$ & \\
      {\bf 13} & (30,4) & $(2, 15, 30)$ &&&&& \\
      \hline
      \hline
      \multicolumn{8}{c}{Genus 7, $\delta = 1$} \\
      \hline
      \hline
      {\bf 14} & (48,48) & $(2, 2, 2, 4)$ & &
      {\bf 15} & (48,41) & $(2, 2, 2, 4)$ & 2 \\
      {\bf 16} & (48,38) & $(2, 2, 2, 4)$ & &
      {\bf 17} & (36,10) & $(2, 2, 2, 6)$ & \\
      {\bf 18} & (32,43) & $(2, 2, 2, 8)$ & &
      {\bf 19} & (32,42) & $(2, 2, 2, 8)$ & 3 \\
      {\bf 20} & (32,39) & $(2, 2, 2, 8)$ & 4 &
      {\bf 21} & (28,3) & $(2, 2, 2, 14)$ & 5 \\
      \hline
      \hline
      \multicolumn{8}{c}{Genus 8, $\delta = 0$} \\
      \hline
      \hline
      {\bf 1} & (336,208) & $(2, 3, 8)$ & &
      {\bf 2} & (336,208) & $(2, 3, 8)$ & \\
      {\bf 3} & (84,7) & $(2, 6, 6)$ & &
      {\bf 4} & (84,7) & $(2, 6, 6)$ & \\
      {\bf 5} & (72,8) & $(2, 4, 18)$ & &
      {\bf 6} & (64,53) & $(2, 4, 32)$ & \\
      {\bf 7} & (60,8) & $(2, 6, 10)$ & &
      {\bf 8} & (48,25) & $(2, 6, 24)$ & \\
      {\bf 9} & (48,17) & $(2, 8, 12)$ & &
      {\bf 10} & (48,28) & $(3, 4, 8)$ & \\
      {\bf 11} & (40,10) & $(2, 10, 20)$ & &
      {\bf 12} & (34,2) & $(2, 17, 34)$ & \\
      \hline
      \hline
      \multicolumn{8}{c}{Genus 8, $\delta = 1$} \\
      \hline
      \hline
      {\bf 13} & (42,1) & $(2, 2, 3, 3)$ & 1, 2, 3, 4 &
      {\bf 14} & (36,4) & $(2, 2, 2, 9)$ & 5 \\
      {\bf 15} & (32,18) & $(2, 2, 2, 16)$ & 6 &
      {\bf 16} & (30,3) & $(2, 2, 3, 5)$ & 7 \\

      \hline
      \hline 
  
      \multicolumn{8}{c}{Genus 9, $\delta = 0$} \\
      \hline
      \hline
      {\bf 1} & (320,1582) & $(2, 4, 5)$ & &
      {\bf 2} & (192,194) & $(2, 3, 12)$ & \\
      {\bf 3} & (192,990) & $(2, 4, 6)$ & &
      {\bf 4} & (192,955) & $(2, 4, 6)$ & \\
      {\bf 5} & (128,138) & $(2, 4, 8)$ & &
      {\bf 6} & (128,136) & $(2, 4, 8)$ & \\
      {\bf 7} & (128,134) & $(2, 4, 8)$ & &
      {\bf 8} & (128,75) & $(2, 4, 8)$ & \\
      {\bf 9} & (120,35) & $(2, 5, 6)$ & &
      {\bf 10} & (120,34) & $(2, 5, 6)$ & \\
      {\bf 11} & (96,187) & $(2, 4, 12)$ & &
      {\bf 12} & (96,186) & $(2, 4, 12)$ & \\
      {\bf 13} & (96,13) & $(2, 4, 12)$ & &
      {\bf 14} & (80,14) & $(2, 4, 20)$ & \\
      {\bf 15} & (72,5) & $(2, 4, 36)$ & &
      {\bf 16} & (64,6) & $(2, 8, 8)$ & \\
      {\bf 17} & (57,1) & $(3, 3, 19)$ & &
      {\bf 18} & (48,5) & $(2, 8, 24)$ & \\
      {\bf 19} & (48,4) & $(2, 8, 24)$ & &
      {\bf 20} & (48,30) & $(4, 4, 6)$ & \\
      {\bf 21} & (42,3) & $(2, 14, 21)$ & &
      {\bf 22} & (40,12) & $(4, 4, 10)$ & \\
      {\bf 23} & (38,2) & $(2, 19, 38)$ & &&&&\\
      \hline
      \hline
    \end{tabular}
      \end{center}
      \end{table}

      \addtocounter{table}{-1}
      \begin{table}
      \caption{(Cont.)}

      \begin{center}
      \renewcommand{\arraystretch}{1.24}
      \begin{tabular}{||c|c|c|c||c|c|c|c||}
      \hline
      \hline
      $\#$ & Group ID & signature & contains & $\#$ & Group ID & signature &
contains \\

\hline
      \hline

      \multicolumn{8}{c}{Genus 9, $\delta = 1$} \\
      \hline
      \hline
      {\bf 24} & (96,193) & $(2, 2, 2, 3)$ & 3 &
      {\bf 25} & (96,227) & $(2, 2, 2, 3)$ & 4 \\
      {\bf 26} & (64,190) & $(2, 2, 2, 4)$ & &
      {\bf 27} & (64,177) & $(2, 2, 2, 4)$ & 5 \\
      {\bf 28} & (64,140) & $(2, 2, 2, 4)$ & &
      {\bf 29} & (64,138) & $(2, 2, 2, 4)$ & 4 \\
      {\bf 30} & (64,135) & $(2, 2, 2, 4)$ & 1, 3, 6 &
      {\bf 31} & (64,134) & $(2, 2, 2, 4)$ & 7 \\
      {\bf 32} & (64,128) & $(2, 2, 2, 4)$ & &
      {\bf 33} & (64,73) & $(2, 2, 2, 4)$ & 2, 8 \\
      {\bf 34} & (48,43) & $(2, 2, 2, 6)$ & 13 &
      {\bf 35} & (48,38) & $(2, 2, 2, 6)$ & \\
      {\bf 36} & (48,15) & $(2, 2, 2, 6)$ & 3 &
      {\bf 37} & (48,48) & $(2, 2, 2, 6)$ & 4, 12 \\
      {\bf 38} & (48,48) & $(2, 2, 2, 6)$ & 11 &
      {\bf 39} & (40,13) & $(2, 2, 2, 10)$ & 14 \\
      {\bf 40} & (40,8) & $(2, 2, 2, 10)$ & &
      {\bf 41} & (36,4) & $(2, 2, 2, 18)$ & 15 \\
      \hline
      \hline

      \multicolumn{8}{c}{Genus 10, $\delta = 0$} \\
      \hline
      \hline
      {\bf 1} & (432,734) & $(2, 3, 8)$ & &
      {\bf 2} & (432,734) & $(2, 3, 8)$ & \\
      {\bf 3} & (360,118) & $(2, 4, 5)$ & &
      {\bf 4} & (324,160) & $(2, 3, 9)$ & \\
      {\bf 5} & (216,92) & $(2, 3, 12)$ & &
      {\bf 6} & (216,158) & $(2, 4, 6)$ & \\
      {\bf 7} & (216,87) & $(2, 4, 6)$ & &
      {\bf 8} & (216,153) & $(3, 3, 4)$ & \\
      {\bf 9} & (180,19) & $(2, 3, 15)$ & &
      {\bf 10} & (168,42) & $(2, 4, 7)$ & \\
      {\bf 11} & (162,14) & $(2, 3, 18)$ & &
      {\bf 12} & (144,122) & $(2, 3, 24)$ & \\
      {\bf 13} & (108,25) & $(2, 6, 6)$ & &
      {\bf 14} & (108,15) & $(2, 4, 12)$ & \\
      {\bf 15} & (88,7) & $(2, 4, 22)$ & &
      {\bf 16} & (80,6) & $(2, 4, 40)$ & \\
      {\bf 17} & (72,28) & $(2, 6, 12)$ & &
      {\bf 18} & (72,23) & $(2, 6, 12)$ & \\
      {\bf 19} & (72,42) & $(3, 4, 6)$ & &
      {\bf 20} & (63,3) & $(3, 3, 21)$ & \\
      {\bf 21} & (60,10) & $(2, 6, 30)$ & &
      {\bf 22} & (42,6) & $(2, 21, 42)$ & \\
      {\bf 23} & (42,2) & $(3, 6, 14)$ & &
      {\bf 24} & (42,2) & $(3, 6, 14)$ & \\
      \hline
      \hline
     
      \multicolumn{8}{c}{Genus 10, $\delta = 1$} \\
      \hline
      \hline
      {\bf 25} & (108,40) & $(2, 2, 2, 3)$ & 4, 6 &
      {\bf 26} & (108,17) & $(2, 2, 2, 3)$ & 1, 2, 7 \\
      {\bf 27} & (72,43) & $(2, 2, 2, 4)$ & 5 &
      {\bf 28} & (72,40) & $(2, 2, 2, 4)$ & 1, 2, 6 \\
      {\bf 29} & (72,15) & $(2, 2, 2, 4)$ & &
      {\bf 30} & (60,5) & $(2, 2, 2, 5)$ & 3, 9 \\
      {\bf 31} & (54,8) & $(2, 2, 2, 6)$ & 11, 14 &
      {\bf 32} & (54,5) & $(2, 2, 3, 3)$ & 5, 8, 13 \\
      {\bf 33} & (48,29) & $(2, 2, 2, 8)$ & 1, 2, 12 &
      {\bf 34} & (44,3) & $(2, 2, 2, 11)$ & 15 \\
      {\bf 35} & (40,6) & $(2, 2, 2, 20)$ & 16 &&&&\\
      \hline
      \hline
      \end{tabular}
      \end{center}
      \end{table}

\clearpage

  \begin{table}[] \label{multicomponents}
      \caption{Reducible Hurwitz loci $\M^*(g,G,\c)$}
      \begin{center}
      \renewcommand{\arraystretch}{1.24}
      \begin{tabular}{||c|c|c|c||}
      \hline
      \hline
       genus & Group ID & signature & components of locus \\
      \hline
      \hline

       7 & (54,6) & (2,6,9) & 6,7 \\
       \hline
       8 & (336,208) & (2,3,8) & 1,2 \\
       \hline
       8 & (84,7) & (2,6,6) & 3,4 \\

       \hline
       9 & (48,48) & (2,2,2,6) & 37,38 \\

       \hline
       10 & (432,734) & (2,3,8) & 1,2 \\
       \hline

       10 & (42,2) & (3,6,14) & 23,24 \\
       \hline
       \hline
      \end{tabular}
      \end{center}
      \end{table}


\bibliographystyle{amsplain}
\bibliography{ref}

\end{document}